\documentclass[11pt,a4paper,oneside,reqno]{amsart}

\usepackage[textwidth=15.5cm]{geometry}

\usepackage{amsmath,amsfonts}
\usepackage{mathtools}
\usepackage[english]{babel}
\usepackage[pdfusetitle]{hyperref}

\newtheorem{theorem}{Theorem}[section]

\newtheorem{corollary}[theorem]{Corollary}
\newtheorem{definition}[theorem]{Definition}

\newcommand{\R}{\mathbb{R}}

\newcommand{\slot}{\,\cdot\,}   
\newcommand{\T}{{\rm T}}        
\newcommand{\D}{{\rm D}}        

\DeclarePairedDelimiter\abs {\lvert}{\rvert}
\DeclarePairedDelimiter\norm{\lVert}{\rVert}

\newcommand{\BC}{C_b}           
\newcommand{\BUC}{C_{b,u}}      

\newcommand{\Ysize}{\eta}       

\newcommand{\sm}{{\scriptscriptstyle M}}
\newcommand{\sx}{{\scriptscriptstyle X}}
\newcommand{\sy}{{\scriptscriptstyle Y}}

\newcommand{\rinj}[1]{r_\text{inj}(#1)}

\title[Persistence of noncompact NHIMs in bounded geometry]%
{Persistence of noncompact normally hyperbolic\\
  invariant manifolds in bounded geometry}

\author{Jaap Eldering}
\email{J.Eldering@uu.nl}
\address{%
  Mathematical Institute,
  Utrecht University,
  Budapestlaan 6,
  3584 CD  Utrecht,
  The Netherlands
}

\date{10 May 2012}

\begin{document}

\maketitle

\begin{abstract}
We prove a persistence result for noncompact normally hyperbolic
invariant manifolds in Riemannian manifolds of bounded geometry.
The bounded geometry of the ambient manifold is a crucial assumption
in order to control the uniformity of all estimates throughout
the proof.
\end{abstract}

\section{Introduction}

Normally hyperbolic invariant manifolds (NHIMs for short) are used in
many areas of dynamical systems, for example, in singular
perturbation theory. It is well-known that compact NHIMs
are persistent under any $C^1$-small perturbation,
see~\cite{Fenichel1971:invarmflds,Hirsch1977:invarmflds}, while
Sakamoto~\cite{Sakamoto1990:invarmlfds-singpert} and Bates, Lu, and
Zeng~\cite{Bates1999:persist-overflow} have extended this to noncompact
NHIMs in Euclidean and Banach spaces respectively. Our result is an
extension to a general noncompact setting in Riemannian manifolds.
Bounded geometry is a crucial additional ingredient, needed to
formulate the necessary uniformity conditions which allow to replace
compactness by uniformity throughout the proof. Bounded geometry can
be viewed as a uniformity condition on the ambient manifold and is
automatically satisfied for Euclidean space.

\section{Bounded geometry}

We follow Eichhorn~\cite{Eichhorn1991:mfld-metrics-noncpt} to define bounded
geometry. Recall that the injectivity radius $\rinj{x}$ at a point
$x \in Q$ is the maximum radius for which the exponential map at $x$
is a diffeomorphism, and that normal coordinates are defined as the inverse map.
\begin{definition}\label{def:bound-geom}
  We say that a complete, finite-dimensional Riemannian manifold
  $(Q,g)$ has $k$-th order bounded geometry when
  \begin{enumerate}
  \item the global injectivity radius $\rinj{Q} = \inf\limits_{x \in Q}\; \rinj{x}$
    is positive, $\rinj{Q} > 0$;
  \item the Riemannian curvature $R$ and its covariant derivatives up
    to $k$-th order are uniformly bounded,
    \begin{equation*}
      \forall\; 0\le i \le k\colon \sup_{x \in Q}\; \norm{\nabla^i R(x)} < \infty,
    \end{equation*}
    with operator norm of $\nabla^i R(x)$ viewed as a multilinear map
    on $\T_x Q$.
  \end{enumerate}
\end{definition}
Note that both Euclidean space and compact smooth Riemannian manifolds
have bounded geometry of any order $k$, i.e.\ $k = \infty$. Less
trivial examples of bounded geometry are symmetric spaces or spaces
constructed as products or as compactly glued connected sums of
bounded geometry spaces.

It follows from Theorem~2.4 in~\cite{Eichhorn1991:mfld-metrics-noncpt} that
a manifold of bounded geometry has an atlas of normal coordinate
charts such that for some fixed $\delta > 0$ there is a normal
coordinate chart defined on each ball $B(x;\delta)$, and moreover, the
representation of the metric $g$ in each chart is $C^k$-bounded,
uniformly over all charts. We shall work with this preferred atlas and
measure the $C^k$ norm of functions in the following way.
\begin{definition}\label{def:unif-bounded-map}
  Let $X,Y$ be Riemannian manifolds of $k\!+\!1$-bounded geometry and
  $f \in C^k(X;Y)$. We say that $f$ is of class $\BC^k$ when there
  exist $\delta_\sx,\,\delta_\sy > 0$ such that for each $x \in X$ we have
  $f(B(x;\delta_\sx)) \subset B(f(x);\delta_\sy)$ and the representation
  \begin{equation}\label{eq:func-coord-repn}
    \tilde{f}_x = \exp_{f(x)}^{-1} \circ f \circ \exp_x
    \colon B(0;\delta_\sx) \subset \T_x X \to \T_y Y
  \end{equation}
  in normal coordinates is of class $\BC^k$ (i.e.\ $C^k$-bounded),
  and moreover, the associated $C^k$-norms of $\tilde{f}_x$ are
  bounded uniformly in~$x \in X$.
\end{definition}
This is a natural definition: $k\!+\!1$-bounded geometry implies that
coordinate transition maps are uniformly $C^k$-bounded, hence this
definition is equivalent to measuring the $C^k$-norm of $f$ at $x$ in
any normal coordinate chart $B(x';\delta_\sx)$ containing $x$. Classes
$\BUC^k(X;Y)$ and $\BUC^{k,\alpha}(X;Y)$ of uniformly (H\"older)
continuous functions can be defined analogously when $X,Y$ are of
$k\!+\!2$-bounded geometry. These ideas can also be extended to
classes $\BC^{k,\alpha}$ and $\BUC^{k,\alpha}$ of vector
fields and submanifolds. We shall allow submanifolds to be
non-injectively immersed.

\section{Results}

We use the following definition of normal hyperbolicity. The flow is
assumed complete for simplicity.
\begin{definition}\label{def:NHIM}
  Let $(Q,g)$ be a smooth Riemannian manifold,
  $\Phi^t \in C^{r \ge 1}$ a flow on $Q$, and let $M \in C^{r \ge 1}$
  be a submanifold of $Q$. Then $M$ is called a normally hyperbolic
  invariant manifold of the dynamical system $(\R,Q,\Phi)$ if all of
  the following conditions hold true:
  \begin{enumerate}
  \item $M$ is invariant, i.e.\ $\forall\; t \in \R\colon \Phi^t(M) = M$;
  \item there exists a continuous splitting
    \begin{equation}\label{eq:NHIM-split}
      \T_M Q = \T M \oplus E^+ \oplus E^-
    \end{equation}
    of the tangent bundle $\T Q$ over $M$ with globally bounded,
    continuous projections $\pi_\sm,\,\pi_+,\,\pi_-$ and this splitting
    is invariant under the linearized flow
    $\D\Phi^t = \D\Phi_\sm^t \oplus \D\Phi_+^t \oplus \D\Phi_-^t$;
  \item there exist real numbers
    $\rho_- < -\rho_\sm \le 0 \le \rho_\sm < \rho_+$ and $C_\sm,C_+,C_- > 0$ such that
    the following exponential growth conditions hold on the various
    subbundles:
    \begin{equation}\label{eq:NHIM-rates}
      \begin{alignedat}{2}
      &\forall\; t\in \R,\,(m,x) \in \T M&&\colon\quad
      \norm{\D\Phi_\sm^t(m)\,x} \le C_\sm\,e^{\rho_\sm\,\abs{t}}\,\norm{x},\\[5pt]
      &\forall\; t\le 0 ,\,(m,x) \in E^+ &&\colon\quad
      \norm{\D\Phi_+^t(m)\,x} \le C_+\,e^{\rho_+\,t}\,\norm{x},\\[5pt]
      &\forall\; t\ge 0 ,\,(m,x) \in E^- &&\colon\quad
      \norm{\D\Phi_-^t(m)\,x} \le C_-\,e^{\rho_-\,t}\,\norm{x}.
      \end{alignedat}
    \end{equation}
\end{enumerate}
\end{definition}

This definition corresponds to `eventual absolute normal
hyperbolicity' in~\cite{Hirsch1977:invarmflds}, and is slightly more
restrictive than the `relative normal hyperbolicity' definition
in~\cite{Hirsch1977:invarmflds} that is also used
in~\cite{Fenichel1971:invarmflds}. We say that $M$ is an
$r$-NHIM if the more general spectral gap condition
\begin{equation}\label{eq:spectral-gap}
  \rho_- < -r\,\rho_\sm \le 0 \le r\,\rho_\sm < \rho_+
  \quad\text{with } r \ge 1
\end{equation}
on the growth exponents above is satisfied.

\begin{theorem}\label{thm:persistNHIM}
  Let $k \ge 2,\, \alpha \in [0,1]$ and $r = k+\alpha$. Let $(Q,g)$ be
  a smooth Riemannian manifold of bounded geometry and
  $v \in \BUC^{k,\alpha}$ a vector field on $Q$. Let
  $M \in \BUC^{k,\alpha}$ be a connected, complete submanifold of $Q$
  that is $r$-normally hyperbolic for the flow defined by $v$, with
  empty unstable bundle, i.e.~$\operatorname{rank}(E^+) = 0$.

  Then for each sufficiently small $\Ysize > 0$ there exists a
  $\delta > 0$ such that for any vector field
  $\tilde{v} \in \BUC^{k,\alpha}$ with
  $\norm{\tilde{v} - v}_1 < \delta$, there is a unique submanifold
  $\tilde{M}$ in the $\Ysize$-neighborhood of $M$, such that
  $\tilde{M}$ is diffeomorphic to $M$ and invariant under the flow
  defined by $\tilde{v}$. Moreover, $\tilde{M}$ is $\BUC^{k,\alpha}$
  and the distance between $\tilde{M}$ and $M$ can be made arbitrarily
  small in $C^{k-1}$-norm by choosing $\norm{\tilde{v} - v}_{k-1}$
  sufficiently small.
\end{theorem}

Let us make some remarks on this result.
\begin{enumerate}
\item%
  The spectral gap condition~\eqref{eq:spectral-gap} of $r$-normal
  hyperbolicity is essential to the proof. The
  $C^{k,\alpha}$ smoothness result is optimal.
  The minimum smoothness requirement $k \ge 2$ is a stronger
  assumption than $k \ge 1$ in the well-known compact case. This seems
  to be intrinsic to the noncompact case, cf.\ hypothesis~H2
  in~\cite{Bates1999:persist-overflow}. If the spectral gap
  condition only holds for some $1 \le r < 2$, then we can
  still obtain a perturbed manifold $\tilde{M}$, but this manifold
  will generally not have better than $C^r$ smoothness.

\item%
  It should be possible to improve this result by lifting some of the
  technical restrictions. First of all, an unstable bundle $E^+$ can
  be added for full normal hyperbolicity. It should hold that the
  persistent manifold $\tilde{M}$ is an $r$-NHIM again.

\item%
  Definition~\ref{def:NHIM} could be relaxed to the more
  general definition of `relative normal hyperbolicity' as used
  in~\cite{Fenichel1971:invarmflds,Hirsch1977:invarmflds,%
    Bates1999:persist-overflow}. This would require using the graph
  transform method; our Perron method proof seems tied to the current
  definition.

\item%
  We only obtain a
  $C^{k-1}$-norm estimate for the perturbation distance of $\tilde{M}$
  away from $M$, even though $\tilde{M} \in C^{k,\alpha}$ is
  preserved. It should be possible to improve this to the perturbation
  being $C^{k,\alpha}$-small when $\norm{\tilde{v} - v}_{k,\alpha}$ is
  small.
\end{enumerate}

By standard phase space extension techniques, we obtain the following
results as a corollary.
\begin{corollary}\label{thm:NHIM-extend}
  Assume the setting of Theorem~\ref{thm:persistNHIM}. If the vector
  field $\tilde{v}$ also depends on time, i.e.
  $\tilde{v} \in \BUC^{k,\alpha}(\R \times Q)$, then there still
  exists a persistent manifold $\tilde{M} \in \BUC^{k,\alpha}$,
  although it may be time-dependent. Similarly, if the vector field
  $\tilde{v}$ depends on an external parameter $p \in R^n$
  and $M$ is an $r$-NHIM for $p = 0$, then there
  exists a neighborhood $U \ni 0$ such that for each $p \in U$ we
  have a unique persistent manifold $\tilde{M}_p \in \BUC^{k,\alpha}$
  and $\tilde{M}_p$ depends $C^{k,\alpha}$ on $p$.
\end{corollary}

\section{Idea of the proof of Theorem~\ref{thm:persistNHIM}}

The following is only a rough sketch of the proof of
Theorem~\ref{thm:persistNHIM}, for a detailed
exposition see~\cite{Eldering2012:persist-noncptNHIM}.

We first reduce the problem to a trivial bundle $X \times Y$ where
$X$ is constructed as a manifold of bounded geometry of sufficiently
high order (say $k+10$) that approximates $M$. Then we embed the
normal bundle $N$ of $X$ into $X \times Y$ with $Y = \R^n$ for some
$n$. A uniform tubular neighborhood of $M$ can be modeled on $N$ since
$M$ is the graph of a small function $h\colon X \to Y$. Additional
normally hyperbolic dynamics can be added in the directions of $Y$
complementary to $N$.

We apply a generalization of the Perron method based on ideas
in~\cite{Henry1981:geom-semilin-parab-PDEs}.
Let $\tilde{v}_\sx$ and $\tilde{v}_\sy(x,y) = A(x)\,y + f(x,y)$ denote
the horizontal and vertical parts of the vector field $\tilde{v}$
respectively. If $(x(t),y(t))$ is a curve in $X \times Y$ with $y(t)$
uniformly small, then we denote by $\Phi_y(t,t_0,x_0)$ the flow of
$\tilde{v}_\sx(\slot,y(t))$ and by $\Psi_x(t,t_0)$ the linear flow of
$A(x(t))$ on $Y$. A contraction map $T$ is defined by
$T(y,x_0) = T_\sy\big(T_\sx(y,x_0),y\big)$ with
\begin{equation}\label{eq:Txy}
  \begin{aligned}
  T_\sx(y,x_0)(t) &= \Phi_y(t,0,x_0),\\
  T_\sy(x,y)(t)   &=
      \int_{-\infty}^t \Psi_x(t,\tau)\,f(x(\tau),y(\tau)) \;{\rm d}\tau
  \end{aligned}
\end{equation}
mappings into appropriate spaces of curves in $X,\,Y$ respectively. We
finally recover the persistent manifold $\tilde{M}$ as the graph of
the map $\tilde{h}\colon x_0 \mapsto \Theta(x_0)(0)$ where $\Theta$
denotes the fixed point of $T$ as a function of the parameter
$x_0 \in X$; the curve $\Theta(x_0)$ in $Y$ is then evaluated at
$t = 0$.

The $\BUC^{k,\alpha}$ smoothness of $\Theta$ is proven inductively
using ideas in~\cite{Vanderbauwhede1987:centermflds} and the fiber
contraction theorem. We introduce certain formal tangent bundles to
work around the problem that the spaces of curves in $X,\,Y$ are not
(Banach) manifolds. We relate holonomy along closed loops in $X$ to
the curvature (which is bounded) to prove uniform continuity of the
formal derivatives of $T$. Restricting the spaces of curves in $X,\,Y$
to bounded time intervals turns these into Banach manifolds; this we
use to finally recover true derivatives that lead to
$\tilde{M} \in \BUC^{k,\alpha}$.

\section*{Acknowledgements}

The author would like to thank Heinz Han\ss mann and Erik van den Ban
for helpful discussions.

\end{document}